\documentclass[12pt,reqno]{amsart}
\usepackage[top=1in,left=1in,bottom=1in,right=1in]{geometry}
\usepackage{amsmath}
\usepackage[square,sort,comma,numbers]{natbib}
\usepackage{amssymb}
\usepackage{amsthm}
\usepackage{tikz}
\usepackage{graphicx}
\usepackage{mathrsfs}
\usepackage{url}
\usepackage{colordvi}
\usepackage{bbm}
\usepackage{color}
\usepackage{verbatim}

\usepackage{xspace}
\usepackage[plain]{fancyref}
\usepackage{placeins}
\usepackage{breqn}
\usepackage{xstring}
\usepackage{enumerate}
\usepackage{natbib}
\usepackage{tkz-graph}
\usepackage[labelfont=normalfont]{caption}
\usepackage{subcaption}

\renewcommand{\k}{\mathbbm{k}}

\renewcommand{\k}{\mathbbm{k}}

\newcommand{\linearspan}[1]{{\langle #1 \rangle}_{\mathbbm{k}}}

\newcommand{\nulla}{\ifmmode{N_{\k}}\else N_{\k}\xspace\fi}

\newtheorem{theorem}{Theorem}[section]
\newtheorem{lemma}[theorem]{Lemma}

\newtheorem{problem}[theorem]{Problem}
\newtheorem{proposition}[theorem]{Proposition}
\newtheorem{corollary}[theorem]{Corollary}
\theoremstyle{definition}
\newtheorem{definition}[theorem]{Definition}
\theoremstyle{remark}
\newtheorem{remark}[theorem]{Remark}
\newtheorem{example}[theorem]{Example}

\numberwithin{subcase}{case}
\numberwithin{subcases}{subcase}

\newcommand*{\fancyrefthmlabelprefix}{thm}
\frefformat{plain}{\fancyrefthmlabelprefix}{theorem\fancyrefdefaultspacing#1}
\Frefformat{plain}{\fancyrefthmlabelprefix}{Theorem\fancyrefdefaultspacing#1}
\newcommand*{\fancyreflemlabelprefix}{lem}
\frefformat{plain}{\fancyreflemlabelprefix}{lemma\fancyrefdefaultspacing#1}
\Frefformat{plain}{\fancyreflemlabelprefix}{Lemma\fancyrefdefaultspacing#1}
\frefformat{plain}{\fancyrefeqlabelprefix}{(#1)}
\Frefformat{plain}{\fancyrefeqlabelprefix}{Equation\fancyrefdefaultspacing#1}
\newcommand*{\fancyrefproplabelprefix}{prop}
\frefformat{plain}{\fancyrefproplabelprefix}{proposition\fancyrefdefaultspacing#1}
\Frefformat{plain}{\fancyrefproplabelprefix}{Proposition\fancyrefdefaultspacing#1}
\newcommand*{\fancyrefdeflabelprefix}{def}
\frefformat{plain}{\fancyrefdeflabelprefix}{definition\fancyrefdefaultspacing#1}
\Frefformat{plain}{\fancyrefdeflabelprefix}{Definition\fancyrefdefaultspacing#1}
\newcommand*{\fancyrefexlabelprefix}{ex}
\frefformat{plain}{\fancyrefexlabelprefix}{example\fancyrefdefaultspacing#1}
\Frefformat{plain}{\fancyrefexlabelprefix}{Example\fancyrefdefaultspacing#1}

\Frefformat{plain}{\fancyreftablabelprefix}{Table\fancyrefdefaultspacing#1}

\makeatletter
\let\@@pmod\pmod
\DeclareRobustCommand{\pmod}{\@ifstar\@pmods\@@pmod}
\def\@pmods#1{\mkern4mu({\operator@font mod}\mkern 6mu#1)}
\makeatother

\tikzstyle{vertex}=[circle, fill, inner sep=0pt, minimum size=5pt]
\newcommand{\vertex}{\node[vertex]}

\begin{document}
\title{Low Degree Nullstellensatz Certificates for 3-Colorability}

\author{Bo Li}
\address{Department of Mathematics.  Harvey Mudd College, Claremont, CA 91711}
\email{bli@g.hmc.edu}

\author{Benjamin Lowenstein}
\address{Department of Mathematics.  Harvey Mudd College, Claremont, CA 91711}
\email{blowenstein@hmc.edu}

\author{Mohamed Omar}
\address{Department of Mathematics.  Harvey Mudd College, Claremont, CA 91711}
\email{omar@g.hmc.edu}


\begin{abstract}
In a seminal paper, De Loera et. al introduce the algorithm NulLA (Nullstellensatz Linear Algebra) and use it to measure the difficulty of determining if a graph is not 3-colorable.  The crux of this relies on a correspondence between 3-colorings of a graph and solutions to a certain system of polynomial equations over a field $\k$.  In this article, we give a new direct combinatorial characterization of graphs that can be determined to be non-3-colorable in the first iteration of this algorithm when $\k=GF(2)$. This greatly simplifies the work of De Loera et. al, as we express the combinatorial characterization directly in terms of the graphs themselves without introducing superfluous directed graphs.  Furthermore, for all graphs on at most $12$ vertices, we determine at which iteration NulLA detects a graph is not 3-colorable when $\k=GF(2)$.
\end{abstract}

\thanks{}
\date{\today}
\maketitle

\section{Introduction}%

In recent years, combinatorial optimization has flourished from algorithms that fundamentally rely on tools from algebraic geometry and commutative algebra.  Work of Lasserre \cite{Lasserre2000}, Lov\'{a}sz-Schrijver \cite{LovaszSchrijver1991}, Sherali-Adams \cite{SheraliAdams1990}, Gouveia, Parrilo and Thomas \cite{GouveiaParriloThomas2008}, and many others have used polynomials to develop approximation algorithms for optimization problems.  Another recent algorithm akin to those above is the Nullstellensatz Linear Algebra algorithm (NulLA) of De Loera et. al \cite{DeLoeraLeeMalkinMargulies08} which addresses feasibility issues in polynomial optimization.  Given a set of polynomials $f_1,f_2,\ldots,f_s \in \k[x_1,\ldots,x_n]$ for some field $\k$, NulLA's goal is to certify that the system of equations $f_1=0, \ f_2=0, \ \ldots \ f_s=0$ has no solution in $\overline{\k}$, the algebraic closure of $\k$.  It exploits Hilbert's Nullstellensatz, a celebrated and fundamental theorem in algebraic geometry (see \cite{CoxLittleOShea92}).

\begin{theorem}[Hilbert's Nullstellensatz]
Let $\k$ be a field, and $f_1,f_2,\ldots,f_s \in \k[x_1,\ldots,x_n]$.  The system of polynomial equations \[f_1=0, \ f_2=0, \ \ldots, \ f_s=0\] has no solution in ${\overline{\k}}^n$ if and only if there are polynomials $\alpha_1,\alpha_2,\ldots,\alpha_s \in \k[x_1,\ldots,x_n]$ such that
\[
1=\alpha_1f_1 + \cdots + \alpha_sf_s.
\]
\end{theorem}

The polynomials $\alpha_1,\alpha_2,\ldots,\alpha_s$ are referred to as a \emph{Nullstellensatz certificate} of infeasibility; indeed they are a witness that the polynomial system $f_1=0, f_2=0, \ldots, f_s=0$ has no solution. The maximum degree of the $\alpha_i's$ is referred to as the \emph{Nullstellensatz degree} of the system, and it is a measure of the complexity of certifying that the system of polynomial equations has no solutions. If a system of polynomial equations is known to have a Nullstellensatz certificate whose Nullstellensatz degree is a small constant (and if $\k$ is finite), one can find a Nullstellensatz certificate in polynomial time in the number of variables through a sequence of linear algebra computations (see \cite{DeLoeraLeeMalkinMargulies08} for details).  However, for general polynomial systems, it is well known that the degree of Nullstellensatz certificates can grow as a function of the number of variables.  

The underlying paradigm in all the above algorithms is the construction of iterative approximations that are tractably computable at early stages.  When applied to combinatorial optimization problems, particularly graph theoretic ones, the key problem that arises is determining the classes of graphs for which a given problem can be resolved in early iterations.  For instance, when applied to the stable set problem, Gouveia, Parrilo and Thomas \cite{GouveiaParriloThomas2008} show that the first iteration of the theta body hierarchy solves the stable set problem for perfect graphs.  This result was first established by Lov\'{a}sz \cite{Lovasz1994} by exploiting polynomials as well.   NulLA itself was originally introduced as a means of unfolding classes of non-3-colorable graphs that can be detected to be non-3-colorable efficiently (that is, in polynomial time in the number of variables of a given graph).  In particular, the authors of \cite{DeLoeraLeeMalkinMargulies08} applied the NulLA algorithm to the following algebraic formulation of graph 3-colorability due to Bayer. We will refer to this as Bayer's formulation throughout the manuscript.

\begin{lemma}[Bayer \cite{Bayer82}] \label{bayer}
  A graph $G$ with vertex set $V$ and edge set $E$ is $3$-colorable if and only if the following system of equations has a solution over an algebraically closed field $\mathbbm{k}$ with $\mbox{char}(\mathbbm{k}) \neq 3$.
\begin{align}
  &0=x_i^3-1& \forall i\in V& \label{eq:vpols} \\
  &0=x_i^2+x_ix_j+x_j^2& \forall v_iv_j \in E&
\end{align}
\end{lemma}

The fundamental concern then is determining combinatorial features of non-3-colorable graphs that dictate the minimum Nullstellensatz degree of infeasibility for the system in Lemma~\ref{bayer}, which we denote by $N_{\k}(G)$, is a small constant.  In light of this, it is natural to address the following problem, a variant of which was first asked in \cite{Omar2010}:

\begin{problem}\label{prob:nulla}
Given a finite field $\k$, and positive integer $d$, characterize those graphs with $N_{\k}(G)=d$.

\end{problem}
Computational evidence (see, for example, Table 1 of \cite{DeLoera2015}) suggests that the minimum Nullstellensatz degree of a non-3-colorable under Bayer's formulation is smallest when the field of coefficients chosen is $GF(2)$ as opposed to $GF(p)$ for primes $p>2$, so from a computational complexity perspective, it may be beneficial to begin addressing Problem~\ref{prob:nulla} by working with Bayer's formulation when $\k=GF(2)$.  A partial answer in this case was given by De Loera et. al \cite{Omar2010} (see their paper for relevant definitions).

\begin{theorem}[Theorem 2.1 of \cite{Omar2010}] \label{thm:origcomb}
For a given simple undirected graph $G$ with vertex set $V=\{v_1,v_2,\ldots,v_n\}$ and edge set $E$, the polynomial system over $GF(2)$ encoding the $3$-colorability of G
\begin{equation*}
  J_G=\{x_i^3+1=0,\ x_i^2+x_ix_j+x_j^2=0:\ i\in V,\ v_iv_j\in E\}
\end{equation*}
has a degree one Nullstellensatz certificate of infeasibility if and only if there exists a set $C$ of oriented partial $3$-cycles and oriented chordless $4$-cycles from $\text{Arcs}(G)$ such that
\begin{enumerate}
  \item $|C_{(v_i,v_j)}|+|C_{(v_j,v_i)}|\equiv0\pmod*2$ for all $v_iv_j\in E$ and
  \item $\sum_{(v_i,v_j)\in \text{Arcs}(G),i<j}|C_{(v_i,v_j)}|\equiv 1 \pmod*2$
\end{enumerate}
where $C_{(v_i,v_j)}$ denotes the set of cycles in $C$ in which the arc $(v_i,v_j)\in\text{Arcs}(G)$ appears.
\end{theorem}

This characterization adds directed structure to undirected graphs, and hence does not fully capture an inherent combinatorial characterization directly from the graphs themselves.  In this paper, we provide such a direct combinatorial characterization for Bayer's formulation when $\k=GF(2)$.  Before introducing the combinatorial characterization, we define the following class of graphs that will play a key role.

\begin{definition}\label{def:deg1graphs}
A graph $G$ with vertex set $V=\{v_1,v_2,\dots, v_n\}$ and edge set $E$ is \emph{covered by length 2 paths} if there exists a set $C$ of length $2$ paths in $G$ such that
\begin{enumerate}
  \item each edge in $E$ appears in an even number of paths in $C$,
  \item the number of paths $v_iv_jv_k$ in $C$ in which $j<i,k$ or $j>i,k$ is odd, and
  \item if $v_i,v_j \in V$ but $v_iv_j \notin E$, then the number of paths in $C$ with $v_i$ and $v_j$ as endpoints is even.
\end{enumerate}
\end{definition}

We now present the main combinatorial characterization.

\begin{theorem} \label{thm:ourcombo}
Let $G$ be a graph. Under Bayer's formulation of 3-colorability with coefficients in the field $\k=GF(2)$, $\nulla(G)=1$ if and only if $G$ is covered by length $2$ paths.
\end{theorem}

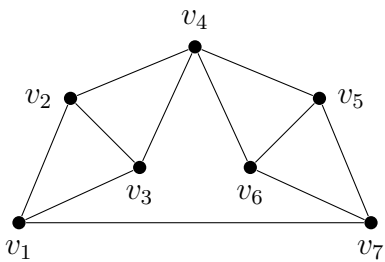
\begin{figure}[ht]
\begin{center}
    \begin{tikzpicture}[scale=1.3]
        \vertex[label=below:$v_1$](v_1) at (180:1.8) {};
        \vertex[label=left:$v_2$](v_2) at (135:1.8) {};
        \vertex[label=below:$v_3$](v_3) at (135:.8) {};
        \vertex[label=above:$v_4$](v_4) at (90:1.8) {};
        \vertex[label=right:$v_5$](v_5) at (45:1.8) {};
        \vertex[label=below:$v_6$](v_6) at (45:.8) {};
        \vertex[label=below:$v_7$](v_7) at (0:1.8) {};
        \path (v_1) edge (v_2) edge (v_3) edge (v_7)
                  (v_2) edge (v_3) edge (v_4)
                  (v_3) edge (v_4)
                  (v_4) edge (v_5) edge (v_6)
                  (v_5) edge (v_6) edge (v_7)
                  (v_6) edge (v_7);
    \end{tikzpicture}
  \caption{Moser spindle.} \label{fig:moser}
\end{center}
\end{figure}

\begin{example}
Let $n$ be a positive integer.  The graph $W_n$, referred to as the \emph{wheel} graph, is the graph whose vertex set is $\{v_1,v_2,\ldots,v_n,w\}$ where the induced subgraph on $V'=\{v_1,v_2,\ldots,v_n\}$ is a cycle and $w$ is a vertex adjacent to all vertices in $V'$.  Without loss of generality, we may assume the cycle whose vertex set is $V'$ has edge set $\{v_1v_2,\ldots,v_{n-1}v_n,v_nv_1\}$.  When $n$ is odd, $W_n$ is not $3$-colorable.  Observe that $C=\{v_1wv_2,v_2wv_3,\ldots,v_{n-1}wv_n,v_nwv_1\}$ is a set of length $2$ paths in $W_n$ that satisfies all the conditions in Definition~\ref{def:deg1graphs}.  We conclude that NulLA (applied to Bayer's formulation when $\k=GF(2)$), detects that odd wheels are non-3-colorable with a degree $1$ Nullstellensatz certificate.  This suggests that NulLA has the potential to be used not only as a computational tool but as a tool for automatic theorem proving.
\end{example}

\begin{example} \label{ex:moserspindle}
Let $G$ be the graph depicted in \Fref{fig:moser}. This non-3-colorable graph is known as the \emph{Moser spindle} and is ubiquitous in graph theory. Using \Fref{thm:ourcombo}, we will show that $\nulla(G)>1$ when $\k=GF(2)$. Suppose otherwise. By \Fref{thm:ourcombo}, $G$ is covered by length 2 paths, which we refer to collectively as $C$. Suppose the path $v_1v_7v_k$ is in $C$ for some vertex $k \notin \{v_1,v_7\}$. Since the edge $v_1v_7$ is not on a $3$-cycle, $v_1v_k$ is not an edge. Moreover, since the edge $v_1v_7$ is not on a $4$-cycle, $v_1v_7v_k$ is the only member of $C$ whose endpoints are $v_1$ and $v_k$. But this is impossible because $v_1v_k$ is not an edge, so the number of paths in $C$ with $v_1$ and $v_k$ as endpoints must be even.  Hence, there are no length $2$ paths in $C$ with $v_1$ and $v_7$ as endpoints.  Thus, through \Fref{thm:ourcombo}, $C$ certifies that $\nulla(G\backslash v_1v_7)=1$. But this is impossible since $G\backslash v_1v_7$ is $3$-colorable, by observation. Thus, no set $C$ with the desired property exists, so $\nulla(G)\neq 1$.
\end{example}

\Fref{ex:moserspindle} generalizes directly in the following way, providing a combinatorial obstruction to existence of a degree 1 Nullstellensatz certificate for Bayer's formulation when $\k=GF(2)$. 

\begin{corollary}\label{cor:corcomb}
Let $\k=GF(2)$, and suppose $G$ is a non-3-colorable graph that contains an edge $e$ for which the following are true:
\begin{itemize}
\item $G\backslash e$ is 3-colorable, and
\item $e$ is not an edge in a $3$-cycle nor a $4$-cycle of $G$.
\end{itemize}
Then $\nulla(G) > 1$.
\end{corollary}

\begin{remark}\label{ex:hard instances}
One of the most celebrated constructions of very hard instances of graph 3-colorability is a construction of Mizuno and Nishihara \cite{Mizuno2008}.  Corollary~\ref{cor:corcomb} is consistent with their findings.  Indeed, in all the graphs they present in Figure 3 (see \cite{Mizuno2008}), the removal of any edge leaves a 3-colorable graph, and each of these graphs has an edge that does not lie on $3$ or $4$-cycle.  This implies that when $\k=GF(2)$, $N_{\k}(G)>1$ for such graphs $G$, so computationally determining that they are not 3-colorable is not immediate under the NulLA paradigm.
\end{remark}

Alongside our combinatorial characterization, in Section 3 we begin the program of determining the Nullstellensatz degree of Bayer's formulation (with coefficients in $GF(2)$) for small non-3-colorable graphs.  Most notably we prove

\begin{theorem} \label{thm:lowdegree}
If $\k=GF(2)$ and $|V(G)| \leq 12$, then $\nulla(G) \leq 4$.
\end{theorem}


\section{Characterizing Degree $1$ Certificates} \label{sec:finitefields}%

This section is dedicated to proving \Fref{thm:ourcombo}, and in particular developing a combinatorial characterization of non-3-colorable graphs $G$ for which $\nulla(G)=1$ (under Bayer's formulation with coefficients in $GF(2)$). We begin with a technical proposition that will be needed throughout:

\begin{proposition} \label{prop:equiv}
  For a graph $G$ with vertex set $V=\{v_1,v_2,\dots,v_n\}$ and edge set $E$, the following are equivalent for Bayer's formulation when $\k=GF(2)$:
\begin{enumerate}
  \item $\nulla(G)=1$.
  \item $1$ is a $\k$-linear combination of
  \begin{align}
    &x_i^3+1& \forall v_i\in V& \\
    &x_k(x_i^2+x_ix_j+x_j^2)& \forall v_iv_j\in E, v_k\in V&
  \end{align}
  \item $1$ is a $\k$-linear combination of
  \begin{align}
    &x_i^2x_j+x_ix_j^2+1& \forall v_iv_j\in E& \\
    &x_i^2x_k+x_ix_jx_k+x_j^2x_k& \forall v_iv_j\in E, v_k\in V, v_i\ne v_k\ne v_j&
  \end{align}
  \item $1$ is a $\k$-linear combination of
  \begin{align}
    &x_i^2x_j+x_ix_j^2+1& \forall v_iv_j\in E& \\
    &x_i^2x_k+x_j^2x_k+x_ix_j^2+x_ix_k^2& \forall v_iv_j\in E, v_jv_k\in E&
  \end{align}
\end{enumerate}
\end{proposition}

The above proposition finds alternate and equivalent sets of polynomials whose solution sets are the same as that of the system in Lemma 1.2.  The last set of polynomials are particularly useful in uncovering our combinatorial characterization. The equivalence of the first three sets was proven in Theorem 2.1 of \cite{Omar2010}. The equivalence to the last set of polynomials follows an argument similar to the proof of Theorem 2.1 in \cite{Omar2010}.  For completeness, we include a proof of this equivalence in the appendix.

In proving \Fref{thm:ourcombo}, we will repeatedly appeal to the following immediate proposition:
\begin{proposition} \label{prop:prop2}

  Let $G$ be a graph with vertex set $V=\{v_1,v_2,\dots, v_n\}$ and edge set $E$, and suppose $G$ is covered by a set $C$ of length $2$ paths. The following statements are equivalent:
\begin{enumerate}
  \item The number of pairs $v_i,v_j\in V$ with $i<j$ for which there are an odd number of paths in $C$ containing $v_i$ with $v_j$ as an endpoint, is itself an odd number.
  \item The sum over all pairs $i<j$ of the number of paths in $C$ containing $v_i$ with $v_j$ as an endpoint is odd.
  \item The number of paths $v_iv_jv_k$ in $C$ in which $j<i,k$ or $j>i,k$ is odd.
\end{enumerate}
\end{proposition}

We now move on to the main result:
\begin{proof}(of \Fref{thm:ourcombo})
  Throughout this proof, for any set of polynomials $S$ in a polynomial ring whose coefficients are in $\k$, we denote by $\langle S \rangle_{\k}$ the linear span of $S$ over $\k$.  Let $F$ be the following set of polynomials:
\begin{align}
  &x_i^2x_j+x_ix_j^2+1& \forall v_iv_j \in E& \label{eq:firstset} \\
  &x_i^2x_k+x_j^2x_k+x_ix_j^2+x_ix_k^2& \forall v_iv_j \in E, v_jv_k \in E& \label{eq:secondset}
\end{align}
By \Fref{prop:equiv}, we know that  $\nulla(G)=1$ if and only if $1\in \langle F\rangle_{\k}$, so we must show that $1\in \langle F\rangle_{\k}$ if and only if $G$ is covered by a set $C$ of length $2$ paths. 

First suppose $G$ is covered by a set of length $2$ paths $C$. Consider the set $H\subset F$ consisting of the following polynomials:
\begin{enumerate}\label{pathtypes}
\item $x_i^2x_k+x_j^2x_k+x_ix_j^2+x_ix_k^2$ for each path $v_iv_jv_k \in C$, and
\item $x_ix_j^2+x_jx_i^2+1$ for each $v_i,v_j\in V$ with $i<j$ such that the number of length 2 paths in $C$ containing $v_i$ and having $v_j$ as an endpoint is odd.
\end{enumerate}


  We claim $1\in\linearspan{H}$ and hence $1\in\linearspan{F}$. Observe that the non-constant monomials appearing in $F$ (and hence in $H$) are all of the form $x_r^2x_s$, where $v_r,v_s \in V$ are arbitrary. We start by showing that the coefficient of $x_r^2x_s$ in $\sum_{h\in H} h$ is $0$ in $\k$, and so all non-constant terms in $\sum_{h\in H} h$ vanish.  An $x_r^2x_s$ term appears in one of four ways:

\begin{enumerate}
  \item[(a)] one $x_r^2x_s$ term for each path in $C$ with $v_r$ and $v_s$ as endpoints, 
  \item[(b)] one $x_r^2x_s$ term for each path in $C$ with $v_r$ as the middle vertex and $v_s$ as an endpoint, 
  \item[(c)] one $x_r^2x_s$ term if there are an odd number of paths in $C$ containing $v_r$ with $v_s$ as an endpoint,
  \item[(d)] one $x_r^2x_s$ term if there are an odd number of paths in $C$ containing $v_s$ with $v_r$ as an endpoint.
\end{enumerate}

The coefficient of the $x_r^2x_s$ in the combined contribution from the terms in (a) and (b) is parity of the number of paths in $C$ containing $v_r$ with $v_s$ as an endpoint.  By our assumption on $C$, there are even number of paths in $C$ containing $v_rv_s$ as an edge, so the coefficient of the $x_r^2x_s$ in the combined contribution from (a) and (b) is the parity of the number of paths in $C$ containing $v_r$ and $v_s$ as endpoints.  The coefficient of $x_r^2x_s$ in the combined contribution from the terms in (c) and (d) is the number of paths in $C$ containing $v_r$ and $v_s$.  By our assumption on $C$, this is again the parity of the number of paths in $C$ containing $v_r$ and $v_s$ as endpoints.  Thus, the coefficient of $x_r^2x_s$ in the combined contribution from the terms in (a),(b),(c),(d) is $0$ in $\k$.  Finally we need to discern the constant term of $\sum_{h \in H} h$.  The constant term is the parity of the number of summands of the form (2).  But this is immediately $1$ by condition (2) of \Fref{def:deg1graphs} and \Fref{prop:prop2}, so the constant term in $\sum_{h \in H} h$ over $\k$ is $1$.


We must now show that if $N_{\k}(G)=1$, then such a set $C$ exists. In that light, Proposition~\ref{prop:equiv} asserts the existence of a set $H$ of polynomials of the form (\ref{eq:firstset}) and (\ref{eq:secondset}) with $\sum_{h \in H} h = 1$. Let $H^*$ be the restriction of $H$ to the polynomials in (\ref{eq:secondset}).  Construct $C$ to consist of the paths $v_iv_jv_k$ for which $x_i^2x_k+x_j^2x_k+x_ix_j^2+x_ix_k^2$ appears in $\sum_{h \in H^*} h$ with a non-zero coefficient.


Suppose $v_rv_s \in E$, define $S_{r,s}$ to be the sum (in $\k$) of the coefficients of the monomials $x_r^2x_s$ and $x_rx_s^2$ appearing in $\sum_{h \in H^*} h$.  Since $\sum_{h \in H} h = 1$, and the only other summand of $\sum_{h \in H} h$ not in $\sum_{h \in H^*} h$ is $x_r^2x_s+x_rx_s^2+1$, $S_{r,s}$ is $0$ in $\k$.  However, the contribution of a single summand $x_i^2x_k+x_j^2x_k+x_ix_j^2+x_ix_k^2$ in $\sum_{h \in H^*}h$ to $S_{r,s}$ is $1$ precisely when $v_rv_s$ is an edge on the path $v_iv_jv_k$, $2$ if $v_r$ and $v_s$ are endpoints of $v_iv_jv_k$, and $0$ otherwise.  Since $S_{r,s}$ is $0$ in $\k$, we deduce that the edge $v_rv_s$ lies on an even number of paths in $C$.

If $v_r,v_s \in V$ but $v_rv_s \not\in E$, then any $x_r^2x_s$ term in $\sum_{h \in H} h$ appears in $\sum_{h \in H^*} h$, so the coefficient of $x_r^2x_s$ in $\sum_{h \in H^*} h$ is $0$ in $\k$. But $x_r^2x_s$ appears once in the summand of $H^*$ corresponding to the path $v_iv_jv_k$ precisely when $v_r$ and $v_s$ are endpoints of $v_iv_jv_k$.  Thus, the number of paths whose endpoints are $v_r$ and $v_s$ is even.

Each edge $v_iv_j$ with $i<j$ contributes a $1$ to the sum $\sum_{h\in H} h$. Moreover, for such pairs $i,j$ the monomial $x_i^2x_j$ appears an odd number of times in $H^*$. We know $x_i^2x_j$ appears in $H^*$ once for each path in $C$ whose endpoints are $v_i$ and $v_j$ and once for each path in $C$ with $v_i$ as a midpoint and $v_j$ as an endpoint. This is equivalent to $x_i^2x_j$ appearing in $H^*$ once for each path in $C$ containing $v_i$ with $v_j$ as an endpoint. Since the number of $1$s appearing in $\sum_{h\in H}h$ is odd, there are an odd number of $i,j$ pairs with $i<j$ such that the number of paths in $C$ containing $i$ with $j$ as an endpoint is odd. By \Fref{prop:prop2}, this establishes condition (2) of \Fref{def:deg1graphs}.

\end{proof}

\section{NulLA on Small Graphs and Future Directions}

Section 2 equipped us with a complete combinatorial understanding of the graphs $G$ for which $N_{\k}(G)=1$ when $\k=GF(2)$.  By Theorem 2.1 of \cite{DeLoera2015}, $N_{\k}(G) \equiv 1 \ (\mbox{mod} \ 3)$, so in investigating Problem~\ref{prob:nulla}, the next natural step is determining when $N_{\k}(G)=4$.  This section is devoted to a systematic study of this for small graphs. We first remark that, in order to find graphs with low minimum Nullstellensatz degree, we only need to focus on a subclass of non-3-colorable graphs.

\begin{definition}[Definition 5.1.4 of \cite{West2001}]
A non-3-colorable graph $G$ is \emph{4-critical} if for any edge $e\in E(G)$, $G\backslash e$ is 3-colorable.
\end{definition}

The following observation is fundamental for our purposes. See Chapter 5 of \cite{West2001} for a discussion of this.

\begin{lemma}[\cite{West2001}] \label{lem:fourcrit}
  Every non-3-colorable graph has a 4-critical subgraph.
\end{lemma}

In light of the previous lemma, the following lemma tells us that 4-critical graphs provide upper bounds for the minimum Nullstellensatz degree of general non-3-colorable graphs. 

\begin{lemma}[Lemma 3.14 of \cite{DeLoeraLeeMarguliesOnn08}] \label{lem:subgraph}
  If $H$ and $G$ are non-3-colorable graphs with $H$ a subgraph of $G$, then $N_{\k}(H)\ge N_{\k}(G)$.
\end{lemma}

\Fref{lem:fourcrit} and \Fref{lem:subgraph} suggest that we solely focus on minimum degree Nullstellensatz certificates for 4-critical graphs. We computed minimum degree Nullstellensatz certificates for all such graphs on at most 12 vertices. A summary of the results are illustrated in \Fref{tab:smallnulla}. With these observations we can now prove \Fref{thm:lowdegree}.

\begin{table}
\begin{center}
\begin{tabular}{c|c|c|c}
$|V|$ & $N_{\k}(G)=1$ & $N_{\k}(G)=4$ & Total $4$-critical graphs \\ \hline
4 & 1 & 0 & 1 \\
5 & 0 & 0 & 0 \\
6 & 1 & 0 & 1 \\
7 & 1 & 1 & 2 \\
8 & 2 & 3 & 5 \\
9 & 5 & 16 & 21 \\
10 & 13 & 137 & 150 \\
11 & 38 & 1183 &1221 \\
12 & 141 & 14440 &14581 \\ \hline
Total & 202 & 15780 & 15982
\end{tabular}
\caption{$N_{\k}(G)$ for $4$-critical graphs on at most $12$ vertices when $\k=GF(2)$.} \label{tab:smallnulla}
\end{center}
\end{table}

\begin{proof}{(of \Fref{thm:lowdegree})}
  By \Fref{lem:fourcrit}, every non-3-colorable graph $G$ on at most 12 vertices has a 4-critical subgraph $H$. By \Fref{tab:smallnulla}, $N_{\k}(H)\in\{1,4\}$. \Fref{lem:subgraph} implies $N_{\k}(G)\leq N_{\k}(H)$. The result then follows by Theorem 2.1 of \cite{DeLoera2015}.
\end{proof}

Many pertinent questions arise from our study of the minimum Nullstellensatz degree of graphs under Bayer's formulation. First and foremost, unless P=NP, one should expect to find a family of graphs for which the minimum Nullstellensatz degree grows arbitrarily large  (see Lemma 3.2 of \cite{DeLoeraLeeMarguliesOnn08} for a discussion of this).  As evidenced by \Fref{tab:smallnulla}, an exhaustive search of the almost 16,000 4-critical graphs on at most 12 vertices indicates that a first step in this direction is to address the following problem:

\begin{problem} \label{thm:highdeg}
For any positive integer $t$, find a graph $G$ for which $N_{\k}(G)>t$.
\end{problem}

After exhaustive experimentation, the authors of \cite{DeLoeraLeeMarguliesOnn08} have yet to see an example resolving Problem~\ref{thm:highdeg} when $t=4$ for any finite field $\k$. One possible method for addressing Problem~\ref{thm:highdeg} is understanding what happens to the minimum Nullstellensatz degree under the famous Haj\'os construction. In his seminal paper \cite{Hajos1961}, Haj\'os defined a recursively constructed class of graphs, which he called $4$-constructible, in the following way:
\begin{enumerate}
\item[i)] $K_4$ is 4-constructible.
\item[ii)] For any two non-adjacent vertices $u$ and $v$ in a 4-constructible graph $G$, the graph obtained from $G$ by adding an edge $e$ incident to $u$ and $v$ and contracting $e$ is also 4-constructible.
\item[iii)] (Haj\'os Construction) For any two 4-constructible graphs $G$ and $H$, with $vw$ an edge of $G$, and $xy$ an edge of $H$, the graph obtained by identifying $v$ and $x$, removing $vw$ and $xy$, and adding the edge $wy$ is also 4-constructible.
\end{enumerate}

Haj\'os proved that the set of 4-constructible graphs is precisely the set of 4-critical graphs, so it is fundamental for us to determine what changes in the minimum Nullstellensatz degree of graphs when applying these constructions. Observe $N_{\k}(K_4)=1$, and by Lemma 3.14 of \cite{DeLoeraLeeMarguliesOnn08}, the minimum Nullstellensatz degree will not increase by applying construction ii). This leads us to the following fundamental question:

\begin{problem}
  Let $G$ and $H$ be 4-critical graphs. What is the relationship between $N_{\k}(G)$, $N_{\k}(H)$ and the minimum Nullstellensatz degree of the graph obtained from $G$ and $H$ by applying the Haj\'os construction?
\end{problem}

\subsection*{Acknowledgments}

The authors thank Jesus De Loera and Susan Margulies for fruitful discussions, the anonymous referees for their helpful feedback, and Eric Stucky for his assistance in typesetting.  

\appendix\section{}

Here, we give a complete proof of Proposition~\ref{prop:equiv}.  The equivalence of the first three sets of polynomials was proven in Theorem 2.1 of \cite{Omar2010}, so we proceed by establishing the equivalence of the fourth set of polynomials with the third.  First, suppose we are given a polynomial $x_i^2x_k+x_j^2x_k+x_ix_j^2+x_ix_k^2$ in (8). Observe that \[x_i^2x_k+x_j^2x_k+x_ix_j^2+x_ix_k^2 = (x_i^2x_k+x_ix_jx_k+x_j^2x_k)+(x_j^2x_i+x_jx_kx_i+x_k^2x_i)\] since char($\k$)=2.  The two summands on the right are polynomials in (6) because $v_iv_j,v_jv_k \in E$,  so any polynomial in (8) is a $\k$-linear combination of polynomials in (6), and hence the fourth condition implies the third.

Now suppose the third condition holds, and we expressed $1$ as a $\k$-linear combination of the polynomials in (5) and (6) as follows:
\begin{equation}\label{longeqn}
1 = \sum_{v_iv_j \in E} \alpha_{ij}(x_i^2x_j+x_ix_j^2+1) + \sum_{\substack{v_iv_j \in E \\ v_k \notin \{v_i,v_j\}, \ k \in V}} \beta_{ijk} (x_i^2x_k+x_ix_jx_k+x_j^2x_k).
\end{equation}
Fix vertices $v_i,v_j,v_k \in V$ with $v_iv_j \in E$.  When the right-hand side of (\ref{longeqn}) is expanded, the coefficient of $x_ix_jx_k$ must be $0$.  We now focus on the contribution of the polynomials in (5) and (6) to coefficient of $x_ix_jx_k$ on the right-hand side for a fixed set of vertices $\{v_i,v_j,v_k\}$.  Without loss of generality, we can demand $v_iv_j \in E$.  

First, suppose neither $v_iv_k$ nor $v_jv_k$ are edges of $G$.  Then the coefficient of $x_ix_jx_k$ is $\beta_{ijk}$.  Comparing both sides of (\ref{longeqn}) implies $\beta_{ijk}=0$, and hence the polynomial $x_i^2x_k+x_ix_jx_k+x_j^2x_k$ does not appear at all on the right-hand side of (\ref{longeqn}).

Now suppose all three of $v_iv_j,v_jv_k,v_iv_k$ are all in $E$.  The polynomials in (5) and (6) that contain $x_ix_jx_k$ are
\[
x_i^2x_k+x_ix_jx_k+x_j^2x_k, \ \ x_k^2x_i+x_ix_jx_k+x_j^2x_i, \ \ x_i^2x_j + x_ix_jx_k + x_k^2x_j,
\]
so if any of these appear as summands of the right-hand side of (\ref{longeqn}), then exactly two of them do (since char($\k$)=2).  If none appear, we do not have to address this case, so assume exactly two appear, and without loss of generality assume they are the latter two.   Then the combined contribution of these summands to the right-hand side of (\ref{longeqn}) is $(x_k^2x_i+x_ix_jx_k+x_j^2x_i )+ (x_i^2x_j + x_ix_jx_k + x_k^2x_j)$.  But observe
\[
(x_k^2x_i+x_ix_jx_k+x_j^2x_i )+ (x_i^2x_j + x_ix_jx_k + x_k^2x_j) = x_i^2x_j+x_k^2x_j + x_ix_j^2+x_ix_k^2
\]
which is a polynomial in (8) since $v_iv_k,v_kv_j \in E$.

Finally, suppose $v_iv_k \notin E, v_jv_k \in E$.  This is the only remaining case since the case when $v_iv_k \in E, v_jv_k \notin E$ follows by symmetry.  In this case, the only polynomials in (5) and (6) that contain $x_ix_jx_k$ as a summand are 
\[
x_i^2x_k+x_ix_jx_k+x_j^2x_k, \ \ x_j^2x_i+x_ix_jx_k+x_k^2x_i.
\]
Again, since the coefficient of $x_ix_jx_k$ must be $0$, either neither of these appear as summands in (\ref{longeqn}) or both do.  Again, we only need consider the case when both appear.  In this case, again since $\k=GF(2)$, $\beta_{ijk}=\beta_{kij}=1$.  Observe then that the contribution of such polynomials to the right hand side of (\ref{longeqn}) is
\[
(x_i^2x_k+x_ix_jx_k+x_j^2x_k )+ (x_j^2x_i+x_ix_jx_k+x_k^2x_i) = x_j^2x_i+x_k^2x_i+x_kx_i^2+x_kx_j^2,
\]
and the latter polynomial is a polynomial in (8) since $v_iv_j,v_jv_k \in E$.  Thus if $1$ is a $\k$-linear combination of polynomials in (6) then it is a $\k$-linear combination of polynomials in (8).

 \bibliographystyle{plain} 
\bibliography{references}

\end{document}